\pdfoutput=1
\documentclass{article}
\usepackage[margin=25mm]{geometry}
\usepackage{amsmath}
\usepackage{amsfonts}
\usepackage{amssymb}
\usepackage{graphicx}

\usepackage[utf8]{inputenc} % правильне кодування
\usepackage[T2A]{fontenc}
\usepackage[english]{babel}
\usepackage[unicode, pdftex]{hyperref}

\usepackage{xcolor}
\usepackage[all]{xy}

\title{Planar graphs as distinguished graph of Morse flows on the 2-disk}
\author{Oleksandr Pryshliak}

\date{\today}

\begin{document}
\maketitle

\begin{abstract}
To investigate the topological structure of Morse flows on the 2-disk we use the planar graphs as destinguished graph of the flow. We assume, that the flow is transversal to the boundary of the 2-disk.
We give a list of all planar graph with at least 3 edges and describe all planar graphs with 4 edges. We use a list of spherical graph with at least 4 edges.
\end{abstract}

\subsection*{Introduction}
The topic of this paper refers to applications of topological graph theory to the classification of Morse flows on surfaces. It uses both topological and algebraic methods when working with embeddings of graphs on the surface.
The main construction, which dates back to the works of Peixoto, consists in the construction of a discriminant graph, the vertices of which are sources, and the edges of which are one-dimensional flowing manifolds.

Classic papers on the topological theory of graphs are books
\cite{angelini2018linear, appel1976every, bang2008digraphs, bollobas1998modern, bondy2008graph, chrobak1994planar, diestel2005graph, dillencourt1993graph, eppstein1999quasi, godsil2001algebraic, gross1999topological, liotta2018beyond, lovasz2001graph, mohar1991embeddings, mohar2001graphs, neumann2015topological, pach2018planar, thomassen2009graphs}.

Embedded  graphs as topological invariants of flows on closed surfaces were constructed in  \cite{bilun2023gradient, Kybalko2018, Oshemkov1998, Peixoto1973, prishlyak1997graphs, prishlyak2020three, akchurin2022three, prishlyak2022topological, prishlyak2017morse, kkp2013, prishlyak2021flows, prishlyak2020topology, prishlyak2019optimal, prishlyak2022Boy},
and on surfaces with a boundary y
\cite{bilun2023discrete, bilun2023typical, loseva2016topology, loseva2022topological, prishlyak2017morse, prishlyak2022topological, prishlyak2003sum, prishlyak2003topological, prishlyak1997graphs, prishlyak2019optimal, stas2023structures}.
For a 3-manifold, a Hegaard diagram is an embedded 4-valent graph in a surface \cite{prish1998vek, prish2001top, Prishlyak2002beh2, prishlyak2002ms, prishlyak2007complete, hatamian2020heegaard, bilun2022morse, bilun2022visualization}.

Morse streams are gradient streams for Morse functions. If we fix the value of the function at special points, then the structure of the flow determines the structure of the function \cite{lychak2009morse, Smale1961}.

Topological invariants, as Reeb graphs, of functions on oriented surfaces were constructed in \cite{Kronrod1950} and \cite{Reeb1946} and in \cite{lychak2009morse} for unoriented surfaces, and in \cite{Bolsinov2004, hladysh2017topology, hladysh2019simple, prishlyak2012topological} for of surfaces with a boundary, in \cite{prishlyak2002morse} for non-compact surfaces.

Embedded graphs as topological invariants of smooth functions were also studied in papers \cite{bilun2023morseRP2, bilun2023morse, hladysh2019simple, hladysh2017topology, prishlyak2002morse, prishlyak2000conjugacy, prishlyak2007classification, lychak2009morse, prishlyak2002ms, prish2015top, prish1998sopr, bilun2002closed, Sharko1993}, for manifolds with a boundary in papers \cite {hladysh2016functions, hladysh2019simple, hladysh2020deformations}, and on 3- and 4-dimensional manifolds in \cite{prishlyak1999equivalence, prishlyak2001conjugacy}.

To get acquainted with the topological theory of functions and dynamical systems, we recommend  \cite{prishlyak2012topological, prish2002theory, prish2004difgeom, prish2006osnovy, prish2015top}.

The purpose of this paper is to describe all possible structures of connected graphs on a plane with no more than 4 edges.

\section{Spherical graphs}
All graphs in this and the next section are connected.
We use spherical graphs to describe all planar graphs. To specify a planar graph, one of the faces must be selected in the spherical graph, which will be external on the planar graph. By projecting a graph from any point of the selected graph by spherical projection onto the plane, we get a flat graph. We will use the list and numbering of spherical graphs given in the work \cite{bilun2023gradient}.

\subsection{Spherical graphs with one edge}

\begin{figure}[ht!]
\center{\includegraphics[width=0.25\linewidth]{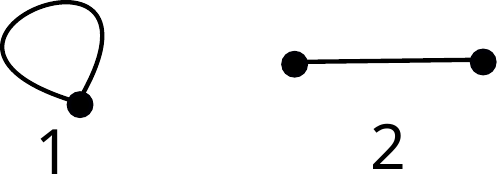}
}
\caption{spherical graphs with one edge}
\label{s1}
\end{figure}

Only two such graphs are possible: a loop (one vertex) and a segment (two vertices) (\ref{s1}).

\subsection{Spherical graphs with two edges}

\begin{figure}[ht!]
\center{\includegraphics[width=0.22\linewidth]{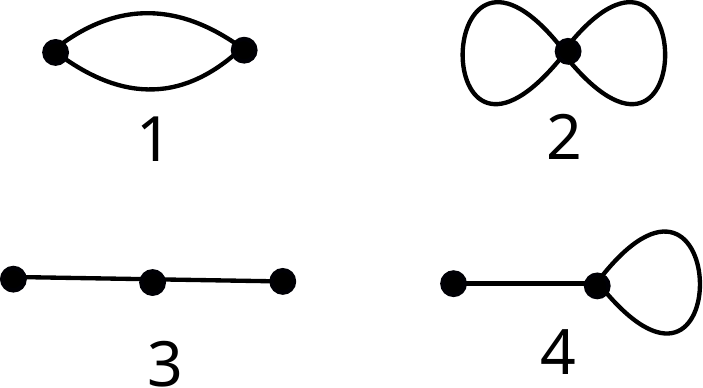}
}
\caption{spherical graphs with two edges}
\label{s2}
\end{figure}

There are 4 spherical graphs with two edges. They are shown in fig. \ref{s2}.

\subsection{Spherical graphs with three edges}

\begin{figure}[ht!]
\center{\includegraphics[width=0.50\linewidth]{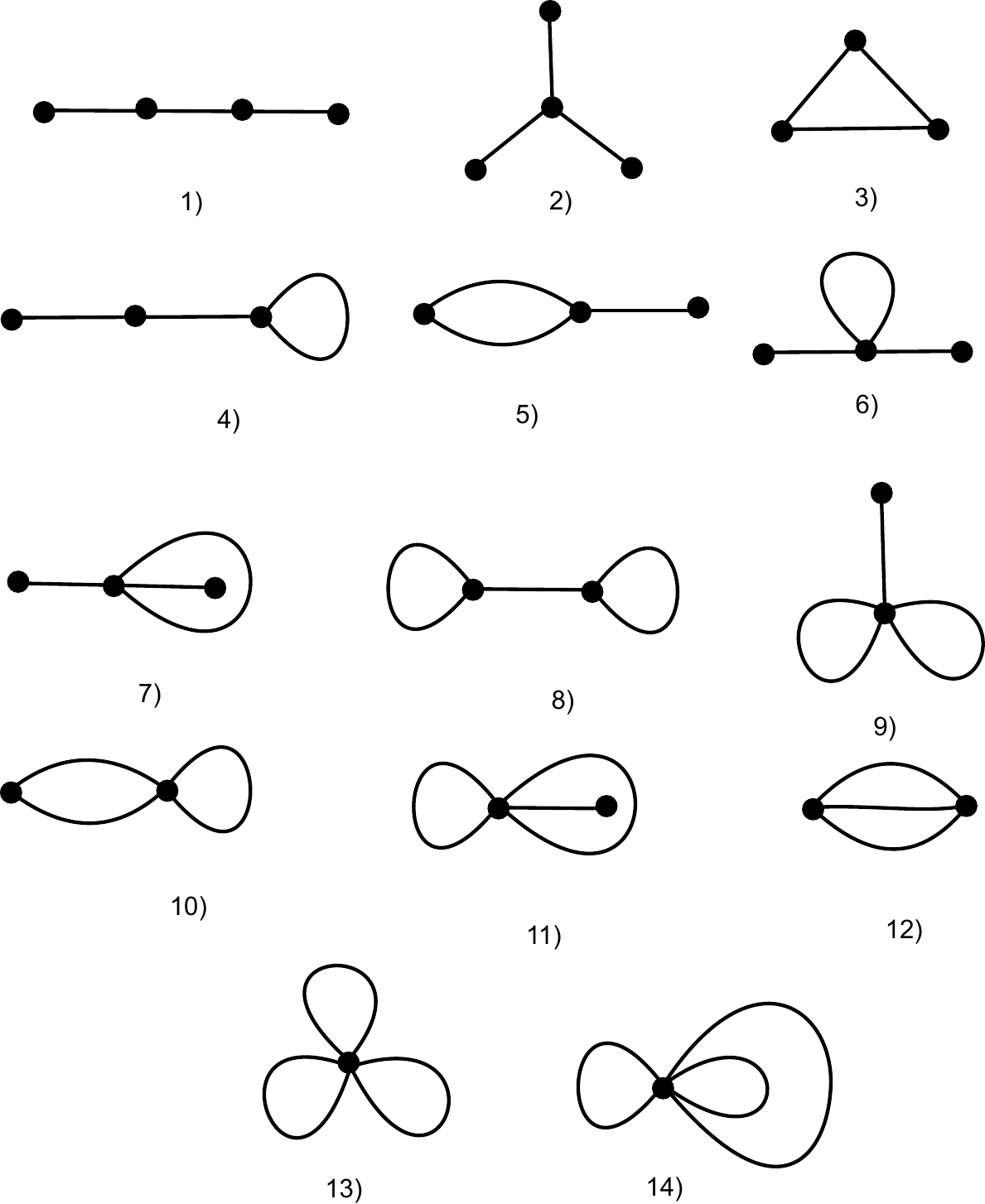}
}
\caption{spherical graphs with three edges
}
\label{s3}
\end{figure}

All spherical graphs with three edges are shown in fig. \ref{s3}.

\subsection{Spherical graphs with 4 edges}
\begin{figure}[ht!]
\center{\includegraphics[width=0.50\linewidth]{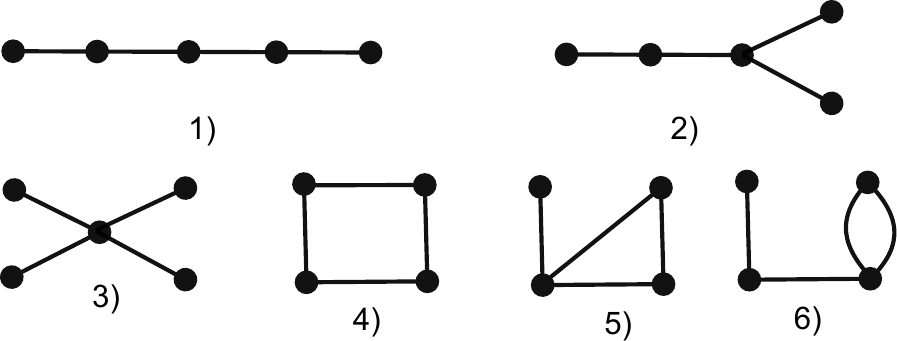}
}
%\caption{spherical graphs with 4 edges (part 1)
%}
%\label{s4a}
%\end{figure}
%\begin{figure}[ht!]
\center{\includegraphics[width=0.50\linewidth]{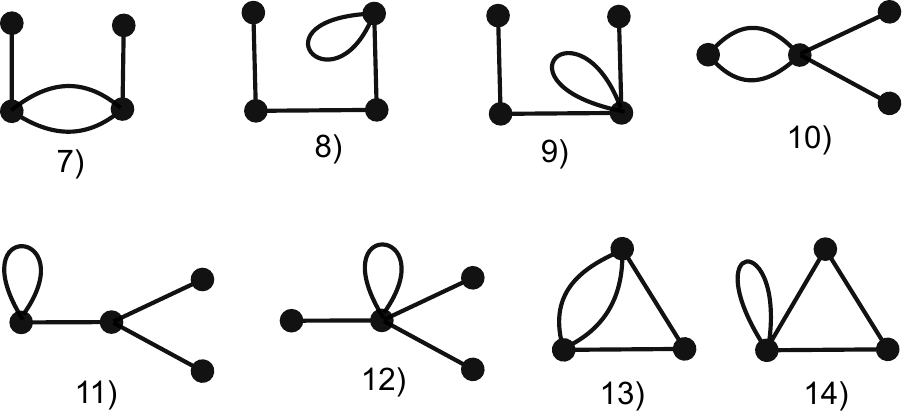}
}
%\caption{spherical graphs with 4 edges (part 2)
%}
%\label{s4b}
%\end{figure}
%\begin{figure}[ht!]
\center{\includegraphics[width=0.50\linewidth]{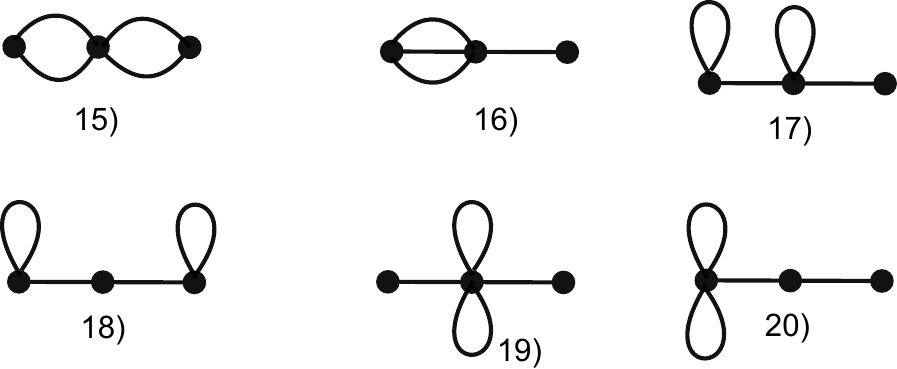}
}
%\caption{spherical graphs with 4 edges (part 3)
%}
%\label{s4c}
%\end{figure}
%\begin{figure}[ht!]
\center{\includegraphics[width=0.50\linewidth]{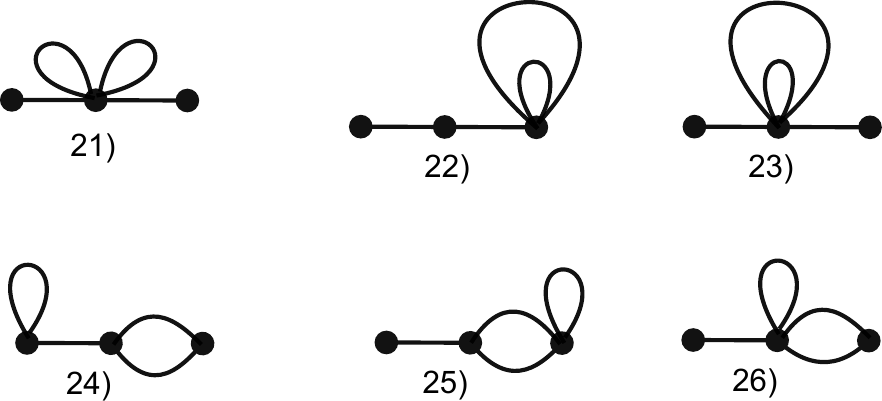}
}
\caption{spherical graphs with 4 edges
}
\label{s4d}
\end{figure}

Spherical graphs with 4 edges are shown in Fig. %\ref{s4a},\ref{s4b},\ref{s4c},
\ref{s4d}.

So, there are a total of 26 non-isomorphic spherical graphs with 4 edges. We took this list and numbering from the work \cite{bilun2023gradient}.

\section{Planar graphs}

The graphs depicted in the previous section are flat (they have a fixed outer edge). The process of changing the outer face can be described as follows: the edge between the layer and the new outer face is torn and connected in a new way, that is, it moves through the point at infinity (the pole of the stereographic projection). At the same time, this process can be repeated several times and with different edges lying within the boundary of the outer face. After obtaining planar graphs in this way, they should be checked for isomorphism among themselves.

\subsection{Planar graphs with one edge}

As for the sphere, only two such graphs are possible: a loop (one vertex) and a segment (two vertices) (\ref{s1}).

If we apply the process of overturning due to infinity to the loop, we get the same graph. This process cannot be applied to a segment.

\subsection{Planar graphs with two edges}

\begin{figure}[ht!]
\center{\includegraphics[width=0.40\linewidth]{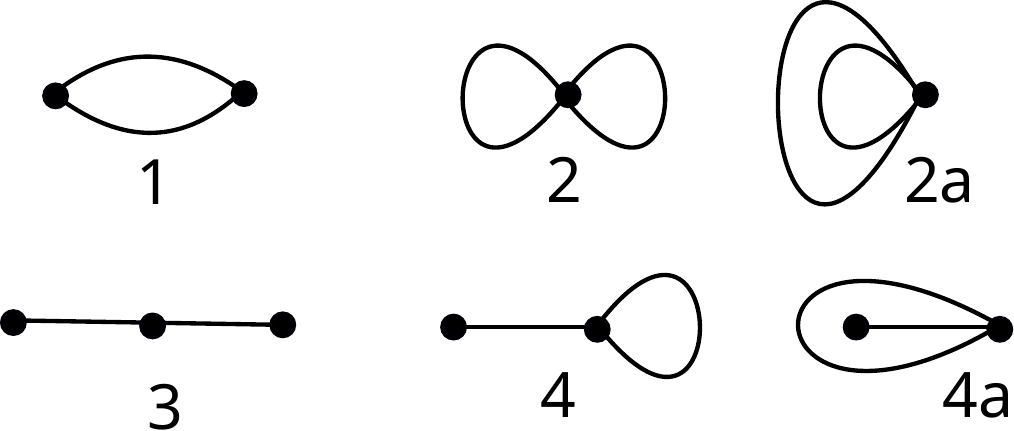}
}
\caption{flat graphs with two edges}
\label{p2}
\end{figure}

There are 6 planar graphs with two edges. They are shown in fig. \ref{p2}. Here, graph 2a is obtained from graph 2 by looping. Similarly, graph 4a is derived from graph 4.

\subsection{Planar graphs with three edges}

\begin{figure}[ht!]
\center{\includegraphics[width=0.65\linewidth]{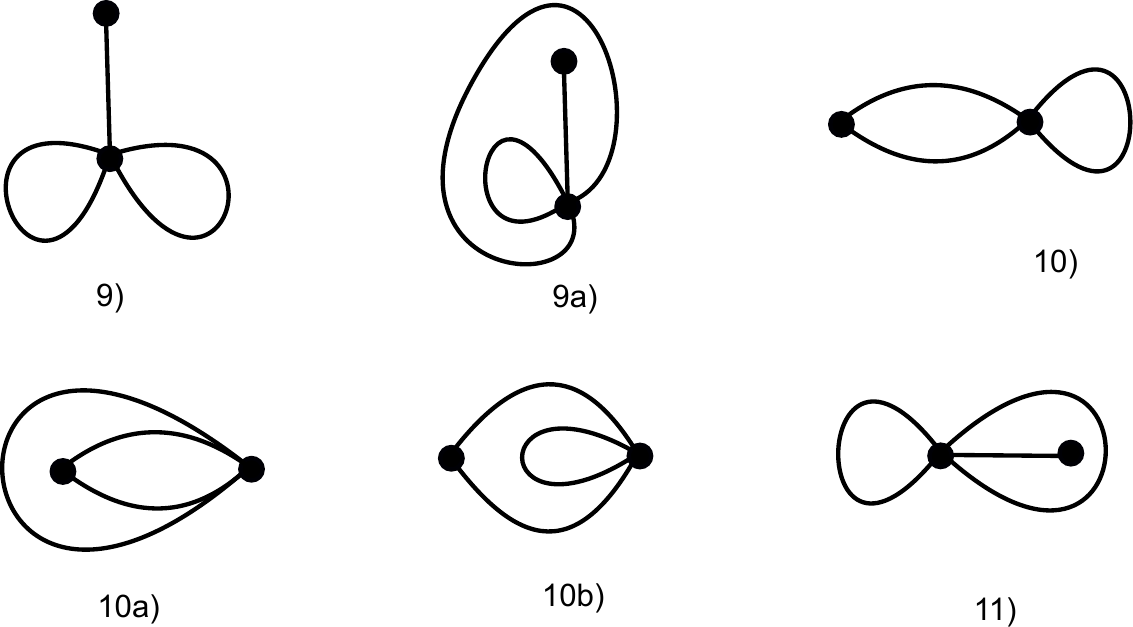}
}
%\caption{spherical graphs with 4 edges (part 3)
%}
%\label{s4c}
%\end{figure}
%\begin{figure}[ht!]
\center{\includegraphics[width=0.65\linewidth]{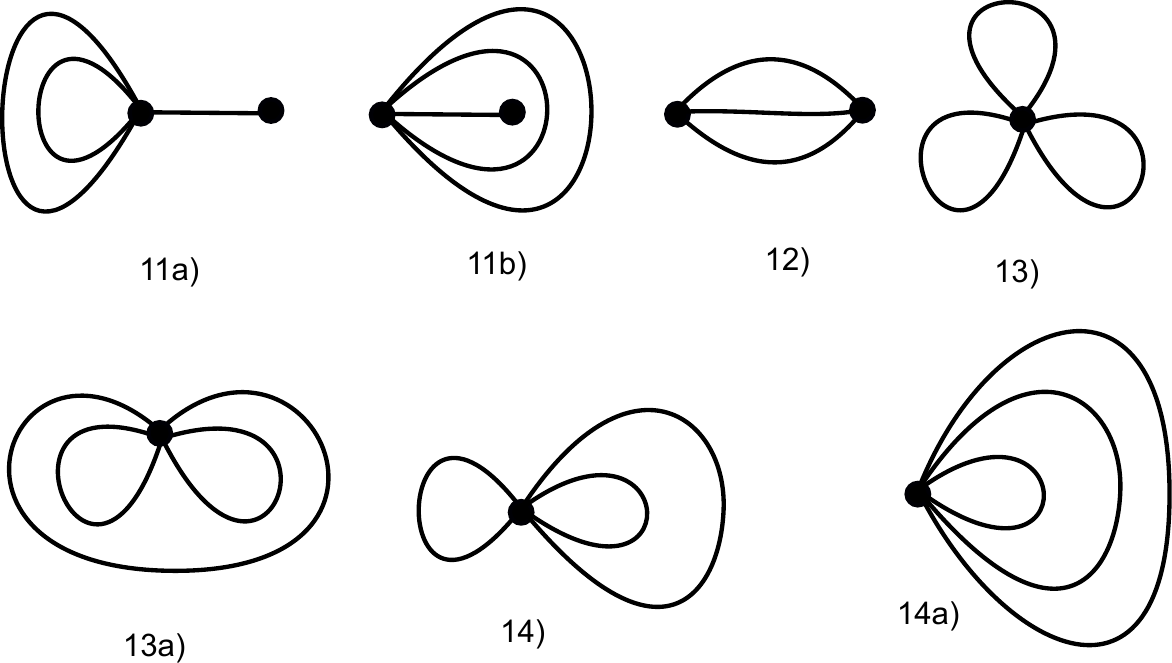}
}
\caption{planar graphs with 3 edges (part 2)
}
\label{p3}
\end{figure}

All planar graphs with three edges are shown in Fig. \ref{p3}. Graphs obtained from each other by flipping edges through infinity have the same numbers but different letters. For example, columns 10a and 10b can be obtained from column 10.

\begin{figure}[ht!]
\center{\includegraphics[width=0.65\linewidth]{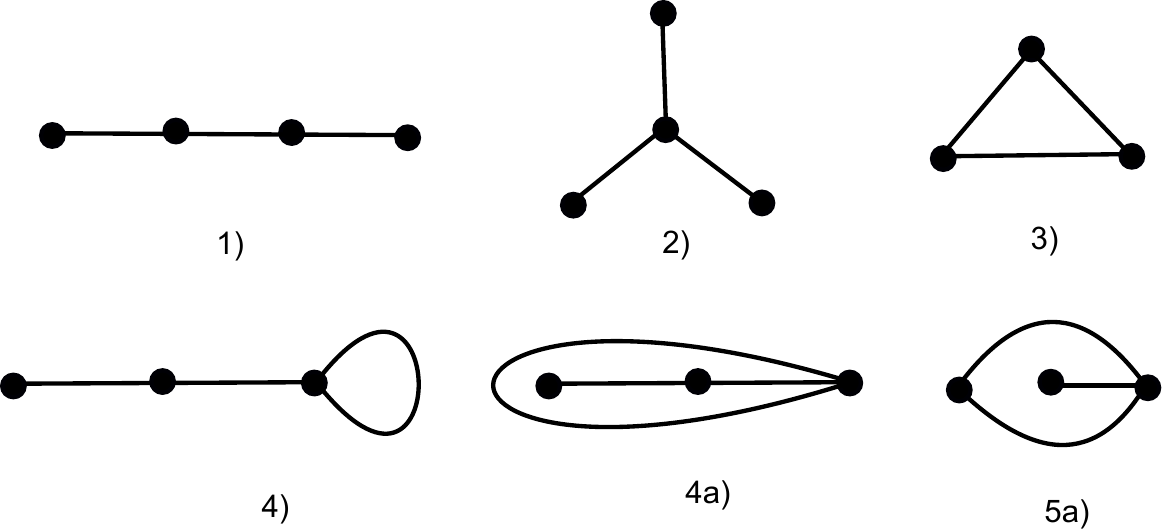}
}
%\caption{spherical graphs with 4 edges (part 1)
%}
%\label{s4a}
%\end{figure}
%\begin{figure}[ht!]
\center{\includegraphics[width=0.65\linewidth]{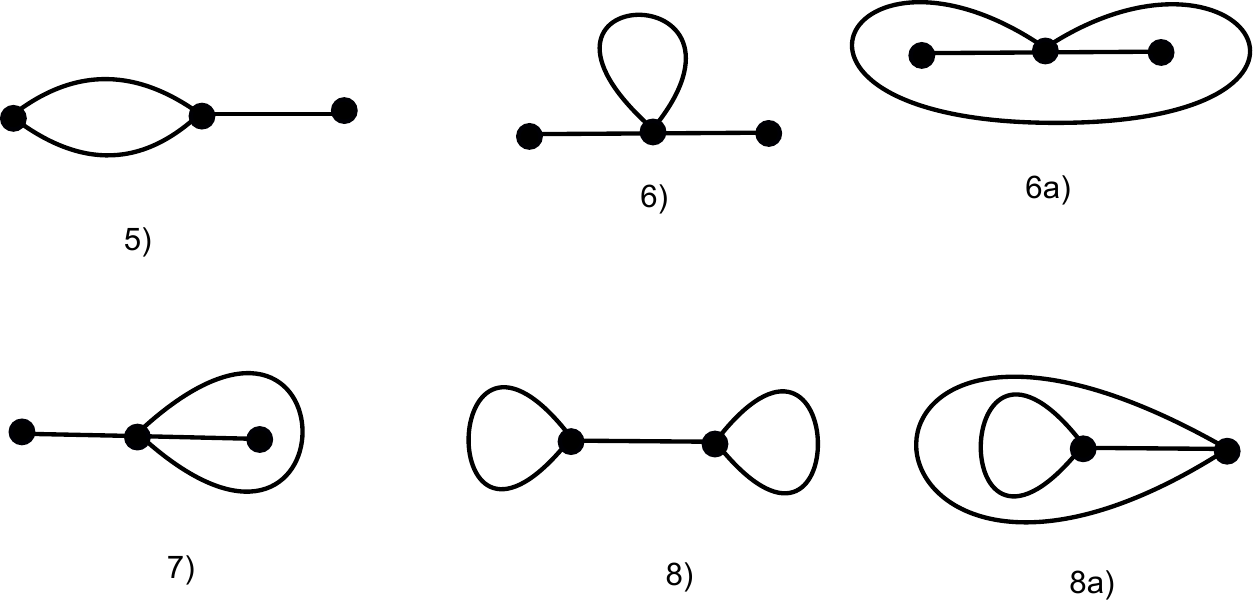}
}
\caption{planar graphs with 3 edges (part 1)
}
\label{p3a}
\end{figure}

So, there are a total of 25 non-isomorphic planar graphs with 3 edges

\subsection{Planar graphs with 4 edges}

To list all planar graphs with 4 edges, we will use the list of spherical graphs with 4 edges shown in Fig. %\ref{s4a},\ref{s4b},\ref{s4c},
\ref{s4d}.

In this figure, the following graphs have a single face (or two faces for which there is an isomorphism that translates one into the other): 1), 2), 3) 4).

Graphs with two asymmetric faces each: 5), 6), 7), 8), 9), 10), 11), 12), 13), 15), 16), 18), 19), 20), 21 ).

Graphs each have three asymmetric faces: 14), 17), 22), 23), 24), 25), 26).

Therefore, there are $4+2\times 15 + 3 \times 7 =55$ non-isomorphic planar graphs with 4 edges.

\section*{Conclusion}
\addcontentsline{toc}{chapter}{Conclusion}

In this paper, we have constructed a complete list of planar graphs with no more than 4 edges. During the construction, we used spherical graphs, fixed the projection face in them, and checked the obtained projections for isomorphism as planar graphs.

The results of the work can be applied in many fields, for example, such graphs are distinguished graphs of gradient vector fields on the 2-disk. So, there two such flow structures with one saddle, 6 structures with two saddles, 25 structures with tree saddles and 55 structures with four saddles.

\newpage
\addcontentsline{toc}{chapter}{Literature}

\newpage
\addcontentsline{toc}{chapter}{Література}

%\bibliographystyle{plain}
%\bibliography{prish}

%\section*{Додатки}
%\addcontentsline{toc}{chapter}{Додатки}
 
%\subsection*{Код даної роботи}
%Ще раз підтверджує те, що кирилиця і код - речі слабосумісні.
%\addcontentsline{toc}{subsection}{Код даної роботи}
 
%\verbatiminput{\jobname.tex} % не має проблем з кирилицею, але має з довгими рядками
%\lstinputlisting{\jobname.tex}
 
\end{document}